\documentclass[12pt]{amsart}

\usepackage{amsmath}
\usepackage{amssymb}
\usepackage{latexsym}
 \newtheorem{theorem}{Theorem}
 
\newtheorem{lemma}[theorem]{Lemma} 
\newtheorem{proposition}[theorem]{Proposition}

\newenvironment{proof*}{\vskip 2mm\noindent {}}{$\Box$ \vskip 2mm}
\numberwithin{equation}{section}

                     \newcommand{\Rr}{{\mathbb{R}}}
                     \newcommand{\Dd}{{\mathbb{D}}}
                     \newcommand{\D}{{\mathbb{D}}}
                     \newcommand{\Zz}{{\mathbb{Z}}}
\renewcommand{\a}{\alpha}
 
\newcommand{\la}{\lambda}

\newcommand{\eit}{e^{i\theta}}

\newcommand{\diz}{(1-|z|^2)}
\newcommand{\diw}{1-|w|^2}

\newcommand{\eps}{\varepsilon}
\renewcommand{\Re}{\mbox{Re}}
\renewcommand{\Im}{\mbox{Im}}

\newcommand{\T}{{\partial \mathbb D}}
             
\thanks{All authors supported by the PICS program no 1019
  of Generalitat de Catalunya
and CNRS. First and second author also supported by European Commission
Research Training Network HPRN-CT-2000-00116.
Second author supported by DGICYT grant
BFM2002-00571 and the CIRIT grant 2001-SGR00431}


\begin{document} 

                     \title{\bf Harmonic and superharmonic majorants on the disk} 

                     \author{Alexander Borichev 
\ \ \ Artur Nicolau \ \ \ Pascal J. Thomas} 

\begin{abstract}
We prove that a positive function on the unit disk admits a 
harmonic majorant if and only if a certain 
logarithmic Lipschitz upper envelope of
it (relevant because of the Harnack inequality) admits 
a superharmonic majorant.  We discuss the logarithmic Lipschitz 
regularity of this superharmonic 
majorant, and show that in general
it cannot be better than that of the Poisson kernel.
We provide examples to show that mere 
superharmonicity of the data does not help with the problem of existence of 
a harmonic majorant.
\end{abstract}

\maketitle

\section{Definitions and statements}

Let $\Dd$ stand for the open unit disk in the complex
plane, and $H^+(\Dd)$ for the cone of positive harmonic
functions on $\Dd$. We would like to describe the
functions $\varphi : \mathbb D \longmapsto \mathbb R_+$
which admit a harmonic majorant, that is,
$h \in H^+ (\mathbb D)$ such that $h \ge \varphi$.
This question arises in problems about the decrease of
bounded holomorphic functions in the unit disk, as well
as in the description of free interpolating sequences for
the Nevanlinna class (see \cite{HMNT}). In that paper,
an answer is given in terms of duality with the measures that act on positive
harmonic functions. The aim of 
this note is to reduce this problem first to the finiteness of a certain 
best Lipschitz majorant function, and then to the existence of a merely 
{\it super}harmonic (nontrivial) majorant.

Let the hyperbolic (or Poincar\'e) distance 
$\rho$ on the disk be defined by
$d\rho(z) :=  (1-|z|^2)^{-1} |dz|$.  This is invariant under biholomorphic 
maps from the disk to itself.  Explicitly, if we first define the 
pseudohyperbolic (or Gleason) distance by
$$
d(z,w):= \left| \frac{z-w}{1-z\overline w} \right|,
$$
then 
\begin{equation}
\label{Poincare}
\rho(z,w)= \frac12 \log \frac{1+d(z,w)}{1-d(z,w)}.
\end{equation}

For $h \in H^+ (\mathbb D)$, the classical Harnack 
inequality reads, for $0<r<1$, $\theta \in \mathbb R$,
$$
\frac{1-r}{1+r} h(0) \le h(r e^{i\theta}) \le \frac{1+r}{1-r} h(0).
$$
This implies  
that the function $\log h$ is Lipschitz with 
constant $2$ with respect to the hyperbolic distance.
We will say that a positive valued function $F$ is \emph{Log-Lipschitz} 
(with constant $C$) if and only if
$|\log F(z) - \log F(w)| \le C \rho(z,w)$ for all $z,w \in \mathbb D$.

\begin{theorem}
\label{superhmaj}
 If  $\varphi$ admits a superharmonic and 
 Log-Lipschitz majorant with constant $C\ge 2$, then 
$\varphi$ admits a harmonic majorant.
\end{theorem}

Since the infimum of two harmonic functions is not in general harmonic, 
there is no smallest harmonic majorant for a given function. On the other 
hand, the cone $Sp(\mathbb D)$ of superharmonic functions is stable under 
finite infima.
Denote by $R(\varphi)$ the reduced function of 
$\varphi$, i.e.
$$
R(\varphi) (z) := \inf \{ u(z) : u \in Sp(\mathbb D), u \ge  \varphi 
\mbox{ on } \mathbb D \}
$$
(see e.g. \cite{Es}). We use the convention that $R(\varphi) \equiv \infty$
if there is no (non identically infinite) superharmonic majorant.
 The reduced function could also be called 
``superharmonic envelope" as in 
\cite{Ra}. 
The reduced function is not in general superharmonic, because it can 
fail to be lower semicontinuous.
J. W. Green 
\cite[Theorem 2]{Gr} proved that when $\varphi$ is continuous, then
$R(\varphi)$, if finite, is also continuous, therefore 
superharmonic, so the infimum in its definition is
really a minimum; furthermore, it is harmonic in the 
open set $\{ z : R(\varphi)(z) > \varphi(z)\}$.

If instead of taking the infimum of all superharmonic functions above 
$\varphi$, we restrict ourselves to those which are Log-Lipschitz with a 
given constant $C$, the corresponding infimum $R_C(\varphi)$
(if finite) will again be Log-Lipschitz with 
constant $C$, hence will be the smallest 
Log-Lipschitz superharmonic majorant of 
$\varphi$ with constant $C$. Theorem \ref{superhmaj} means that $\varphi$ 
amits a harmonic majorant if and only if $R_C(\varphi)$ is finite. 
 
In order to study the finiteness 
of $R_C(\varphi)$, it would be nice to be able
 to proceed in two steps, first dealing 
with the Lipschitz property, then with  superharmonicity. 
The smallest Log-Lipschitz majorant with constant $C$ of a given 
nonnegative
function $\varphi$ is (see \cite{McS}, or \cite{JBHU} 
for a more recent survey)
\[
L_C(\varphi) (z) 
= \exp \{ \sup_{w\in \mathbb D} (\varphi(w)- C \rho(w,z)) \}.
\]
Clearly, $L_C(\varphi) \le R_C(\varphi)$, thus
$R(L_C(\varphi)) \le R_C(\varphi)$. 
That inequality can be strict when $C<2$.

\begin{proposition} 
\label{RnotLip}
For every $0<\gamma\le\beta<2$ there exists
a bounded positive function $H$ on the disc such that
$$
R(L_\gamma(H))\ne R_\beta(H).
$$
\end{proposition}

But we don't know
whether $R(L_C(\varphi))=R_C(\varphi)$ holds for $C \ge 2$. 
Nevertheless, if a Log-Lipschitz function admits a superharmonic majorant 
$s$, then the invariant averages of $s$ will provide us with a 
superharmonic majorant with a weak Log-Lipschitz property, and 
finally we obtain the following 
result.

\begin{theorem}
\label{twostepmaj}
If there exists $C >0$ such that $L_C (\varphi)$ admits a 
superharmonic majorant, then
$\varphi$ admits a harmonic majorant.
\end{theorem}

Note that it would not help to perform our two steps in the reverse order. 
Typical data for many problem of harmonic majorants are functions 
$\varphi$ which 
vanish everywhere except on a discrete subset \cite{HMNT}. For such $\varphi$, 
$R(\varphi) = \varphi$.

\vskip.5cm
Of course, Theorem \ref{twostepmaj} 
begs the question of criteria to ensure the finiteness of 
the reduced function of Log-Lipschitz data. In this general 
direction, one should note results of Koosis \cite[p. 77]{Ko} and Cole and 
Ransford \cite[Theorem 1.3]{CR} which imply that 
when $\varphi$ is merely continuous,
$$
R\varphi (x) = \sup\{ \int \varphi d\mu : \mu \in I_x \}
=  \sup\{ \int \varphi d\mu : \mu \in H_x \}, 
$$
where $I_x$ denotes the set of Jensen measures for $x$, 
and $H_x$ denotes the set of harmonic measures for  $x$ 
with respect to a domain $\omega \subset \subset \Dd$.
A perhaps more computable characterization of $R\varphi$ is given by 
\cite[Th\'eor\`eme, p. 80]{Ko} : let $D_H(z,r)$ stand for the disc of 
center $z$ and radius $r$ with respect to the hyperbolic distance $\rho$,
and $d\beta (z) := \diz^{-2} d \la_2 (z)$ be the invariant 
measure on the disk. 
Given a real-valued continuous function on $\Dd$, let
\[
MF(z) := \sup_{r>0} \frac1{\beta (D_H(z,r))}\int_{D_H(z,r)} F(w) d\beta(w),
\]
and define $F^{(0)}:= F$, $F^{(k+1)}:= MF^{(k)}$. Then, arguing  as 
in \cite{Ko}, one can check that
\[
R \varphi (z) = \lim_{n\to\infty} \varphi^{(n)}(z).
\]
So, Theorem \ref{twostepmaj} says that $\varphi$ admits a harmonic 
majorant if and only if the sequence $L_C (\varphi)^{(n)}(0)$ remains 
bounded.  

\vskip.5cm

{\bf Counterexamples.}

We want to see that the a priori assumption of superharmonicity upon the data, 
in contrast to the assumption of a Lipschitz type property,
does not improve the prospect for the existence of a harmonic majorant. 
A result of that type can be given as follows. First recall that a 
necessary condition for the existence of a harmonic majorant can be given 
in terms of nontangential maximal functions. 

The Stolz angle with vertex at some point $\zeta \in \partial\Dd$ is 
defined by 
$$
\Gamma_\alpha (\zeta) := \{ z \in \mathbb D : |z-\zeta|
\le \alpha (1-|z|^2) \}.
$$
Given a function $f$ from $\Dd$ to $\mathbb R_+$,
the nontangential maximal function is defined as
\[
M f (\zeta) := \max_{\Gamma_\alpha (\zeta)} f .
\]
Finally, for a function $g$ from some measure space $(X,\mu)$ to  $\mathbb R_+$,
we write that $g \in L^1_w (X)$ (``weakly integrable functions'') if and only if 
there is $C>0$ such that 
\[
\mu(g \ge \lambda) \le C/\lambda, \mbox{ for all }
\lambda >0. 
\]
It is a well-known and easy fact that if $\varphi$ admits a 
harmonic majorant, then $M\varphi \in L^1_w(\partial \mathbb D)$ (see 
\cite[Proposition 1.5]{HMNT}, \cite[Lemma 14, p. 6]{PT1}, 
\cite[Lemma 5, p. 6]{PT2}). This condition 
can't be improved in a quantitative way, even for superharmonic data.

For notational convenience, the examples below are given in the upper-half 
plane $U_+ := \{ z \in \mathbb C : \Im z >0\}$.  Stolz angles are then 
defined as
$$
\Gamma_\alpha (t) := \{ z = x+iy \in U_+ : |x-t|
\le \alpha y \}.
$$

\begin{proposition}
\label{sharpmaxf}
For any decreasing function $s$ from $\Rr_+^*$ to $\Rr_+^*$ such that 
$\lim_{t\to\infty} s(t) =0$ and
$\int_1^\infty s(t) dt = \infty$, there exists a positive superharmonic 
function $\varphi$ on $U_+$ such that 
\[
\left| \{ \xi \in \Rr : M \varphi (\xi) > t \} \right| \le s(t)
\mbox{ for any } t>0 \mbox{ large enough,}
\]
but there does not exist any harmonic function $h\in Ha (U_+)$ such that
$h(z) \ge \varphi (z)$ for all $z \in U_+$.
\end{proposition}

Here, as in various further instances,
 $|\cdot|$ stands for one-dimensional
Lebesgue measure on the real line or on the circle.

Examples of the same type allow us to show that no condition relying 
solely on the rate of increase of $\varphi$ near the boundary can be 
sufficient to imply that $\varphi$ is dominated by some harmonic function,
even assuming as above that $\varphi$ be superharmonic.

\begin{proposition}
\label{anyrate}
For any decreasing function $s$ from $\Rr_+^*$ to $\Rr_+^*$ such that 
$\lim_{t\to 0} s(t)= \infty$, there exists a positive superharmonic 
function $\varphi$ on $U_+$ such that 
\[
\varphi(x+iy) \le s(y), \mbox{ for any } x+iy \in U_+ , y \mbox{ small 
enough },
\]
but there does not exist any harmonic function $h\in Ha (U_+)$ such that
$h(z) \ge \varphi (z)$ for all $z \in U_+$.
\end{proposition}

Observe that Proposition \ref{anyrate} would be very 
easy if we allowed $\varphi$ to be 
subharmonic; then it would be enough to select $\varphi(x+iy) := s_1(y)$, with 
$s_1$ the lower convex envelope of $s$. This remark is made to stress the 
difference between our problem and the more classical question of harmonic 
majorants of subharmonic functions.

\vskip.5cm
The paper is organized as follows. The next section is devoted to studying 
the dyadic analogue of the problem of harmonic majorants. Theorems \ref{superhmaj} and 
\ref{twostepmaj} are proved in sections \ref{shm} and \ref{compenv} respectively.
Proposition \ref{RnotLip} is proved in section \ref{PRNL}. The other
counterexamples are given in the last section of the paper.

\section{A discrete model}

Recall that any positive harmonic function $h$ is the Poisson integral of a 
finite positive measure $\mu$ on the boundary of the disk,
$h(z) = \int_0^{2\pi} P_z (\eit) \, d\mu(\theta)$, where
\[
 P_z (\eit) := \frac1{2\pi} \frac{\diz}{|z-\eit|^2}.
\]
The following considerations concern the simpler case of functions which 
are generated by the ``square" kernel
\[
K_z (\eit)  :=\frac1{|I_z|} \chi_{I_z}(\eit),\mbox{ where }
I_z := \{ \eit : z \in \Gamma_\alpha (\eit) \}.
\]
Here  $\chi_E$ stands for 
the characteristic function of the set $E$. The ``square" integral of a 
finite measure $\mu$ is defined by 
$\int_0^{2\pi} K_z (\eit) \, d\mu(\theta) = \mu(I_z)/|I_z|$.

Consider the usual partition of $\partial \Dd$ in dyadic arcs,
for any $n$ in $\Zz_+$:
\[
I_{n,k} := \{ \eit : \theta \in [{2\pi}k 2^{-n},{2\pi}(k+1) 2^{-n})
\}, \ 0 \le k < 2^n .
\]
Note that $|I_{n,k}|= 2 \pi 2^{-n}$.
To this subdivision we associate the Whitney partition in ``dyadic squares" 
of the unit disk :
\[
Q_{n,k} := \{ r \eit : \eit \in I_{n,k} , 1- 2^{-n} \le r < 1- 2^{-n-1} \}.
\]
It is well known and easy to see that there exists a constant $c_\a$ such 
that for any $z \in Q_{n,k}$,$P_z \ge c_\a K_{I_{n,k}}$.
This implies that sufficient conditions for majorization by ``$K$-harmonic 
functions" yield sufficient conditions for majorization by (true) harmonic 
functions.

\begin{theorem}
\label{discrete}
Given a collection of nonnegative data $\{p_{n,k}\} \subset \Rr_+$, there 
exists a finite positive measure $\mu$ on $\partial \Dd$ such that 
\[
\frac{ \mu(I_{n,k})}{|I_{n,k}|} \ge p_{n,k}
\]
if and only if there exists a constant $S$ such that
\begin{equation}
\label{discNSC}
\sum p_{n,k} |I_{n,k}| \le S 
\end{equation}
where the sum is taken over any disjoint 
family of dyadic arcs $\{I_{n,k}\}$.
\end{theorem}

\begin{proof}
The condition is clearly necessary with $S=\mu(\partial\Dd)$.
To prove the converse direction, let us consider the following modified 
data:
\[
\tilde p_{n,k} := \frac1{|I_{n,k}|} \sup
\left\{ \sum p_{q,j} |I_{q,j}| 
\right\} ,
\]
where the supremum is taken over any disjoint 
family  $\{I_{q,j}\}$ of dyadic subintervals of $ I_{n,k}$. Observe that 
$\tilde p_{n,k} \ge p_{n,k}$, and that $\tilde p_{n,k}$ verifies the 
following discrete superharmonicity property :
\[
\tilde p_{n,k} \ge \frac12 ( \tilde p_{n+1,2k} + \tilde p_{n+1,2k+1} ).
\]

Assuming (\ref{discNSC}) we will construct 
a sequence of positive measures $\mu_n$ of bounded total mass,
uniformly distributed on each interval $I_{n,j}$, such that 
for all $m \le n$,
\begin{equation}
\label{hyprecmes}
\mu_n(I_{m,k}) \ge |I_{m,k}| \tilde p_{m,k}.
\end{equation}

Let $\mu_0$ be the uniform measure of total mass $S$ on the interval 
$I_{0,0}= \partial \Dd$.  The hypothesis (\ref{discNSC}) coincides with
(\ref{hyprecmes}) in this case.
Suppose that $\mu_m$, $m \le n$ have already 
been constructed satisfying (\ref{hyprecmes}), we will choose $\mu_{n+1}$. 
Fix $j$, $0 \le j < 2^n$. Then $I_{n,j} = I_{n+1,2j} \cup I_{n+1,2j+1}$, 
so (\ref{hyprecmes}) implies in particular that
\[
\mu_n (I_{n,j} ) \ge |I_{n+1,2j}| \tilde p_{n+1,2j} +  |I_{n+1,2j+1}| 
\tilde p_{n+1,2j+1},
\]
so we find $\a$, $\beta \ge 0$ such that 
\[
\mu_n (I_{n,j}) = \a + \beta, \quad
\a \ge |I_{n+1,2j}| \tilde p_{n+1,2j}, \quad 
\beta \ge |I_{n+1,2j+1}| \tilde p_{n+1,2j+1}.
\]
We set $\mu_{n+1} (I_{n+1,2j}) = \a$, $\mu_{n+1} (I_{n+1,2j+1}) = \beta$. 
This defines a measure $\mu_{n+1}$ which verifies (\ref{hyprecmes}) at 
rank $n+1$ for $m=n+1$. It also
satisfies 
\[
\mu_{n+1} (I_{m,j}) = \mu_{n} (I_{m,j}), \quad \forall m \le n,
\]
so (\ref{hyprecmes}) is verified by 
$\mu_{n+1}$ for all $m \le n+1$.  
This bounded sequence of measures contains a weakly convergent 
subsequence, whose limit $\mu$ will satisfy the condition in the 
Theorem.
\end{proof}

\section{Proof of Theorem \ref{superhmaj}}
\label{shm}

We will prove the following slightly more general fact, which 
will be useful in Section \ref{compenv}. 

\begin{proposition}
\label{liplarghmaj}
Let $u$ be a positive superharmonic function on the disk satisfying that
for any $z$, $w \in \Dd$,
\begin{equation}
\label{LipinLarg}
\left| \log v(z) - \log v(w) \right| \le C_1 (1+ \rho(z,w)).
\end{equation}
Then there exists $h \in H^+ (\Dd)$ such that
$h \ge u$.
\end{proposition}

\begin{proof}

The Riesz representation theorem tells us that there exists a positive measure 
$\nu = - \Delta u$ in the sense of distributions, and a positive harmonic function 
$h_0$ such that
\[
u(z) = h_0(z) + \int_{\Dd} \log \frac1{d(z,w)} d \nu(w) =: 
 h_0(z) + u_0(z).
\]
Let $\delta \in (0,1)$ to be chosen later, then
\begin{multline*}
u_0(z) = \int_{\{w: d(z,w) \le \delta\}} \log \frac1{d(z,w)} d \nu(w)
+  \int_{\{w: d(z,w) > \delta\}} \log \frac1{d(z,w)} d \nu(w)
\\
=: u_1(z) + u_2(z) .
\end{multline*}
For $d(z,w) > \delta$, we have
\[
\log \frac1{d(z,w)} \le C_\delta \frac{1-|zw|^2}{|1-z\overline w|^2} 
(\diw) ,
\]
which is harmonic in $z$. Note that 
\begin{equation}
\label{Blaschkemeas}
\infty > u(0) \ge \int_{\Dd} \log \frac1{|w|} d \nu(w)
\ge \int_{\Dd} (1-|w|) d \nu(w),
\end{equation}
so that the harmonic fuction
\[
h_2(z) :=  
 C_\delta \int_{\Dd} \frac{1-|zw|^2}{|1-z\overline w|^2} (1-|w|^2) d \nu(w)
\]
is bigger than $u_2$. Now we only need to find a harmonic majorant for the 
remaining term $u_1$.

A sequence $\{z_k\} \subset \Dd$ is called \emph{uniformly dense}
if there 
exists $0<r_1<r_2$ satisfying
\begin{enumerate}
\item $D_H (z_j,r_1) \cap D_H (z_k,r_1) = \emptyset$ for any $j\neq k$,
\item $\Dd \subset \cup_k D_H (z_k,r_2)$.
\end{enumerate}
Recall that $D_H (z,r)$ is the hyperbolic disk with center $z$ and radius $r$.

\begin{lemma}
\label{majnet}
For $\delta \in (0,1)$ a properly chosen absolute constant, 
there exist a uniformly dense
sequence $\{z_k\}$ and a positive harmonic function $h_1$ such 
that for any $k$, $u_1(z_k) \le h_1(z_k)$.
\end{lemma}

Accepting this Lemma, we can easily finish the proof : if we write $h_3 := 
h_0+h_1 + h_2$, we see that we have $u(z_k) \le h_3(z_k)$, for any $k$. 
Now by Harnack's inequality, for any $z \in D_H (z_k,r_2)$,
$h_3(z)\ge e^{-2r_2} h_3(z_k)$;
while by (\ref{LipinLarg}) ---note that this is the only step in this 
argument where this hypothesis is used--- 
\[
u(z) \le \exp[C_1(1 + r_2)] u(z_k) \le \exp[C_1(1 + r_2)] h_3(z_k) 
\le \exp[C_1(1 + r_2)] e^{2r_2} h_3(z) .
\]
So we have found a harmonic majorant of $u$.
\end{proof}

\begin{proof*} {\bf Proof of Lemma \ref{majnet}}

First let $T_{n,j} := \left\{ z \in Q_{n,j} : \rho(z, \partial  Q_{n,j} ) > \delta 
\right\}$.
For $\delta$ small enough, 
denoting by $m$ the $2$-dimensional Lebesgue measure on the plane, 
\[
 m (Q_{n,j}) \ge m (T_{n,j}) 
\ge C(\delta) m (Q_{n,j}) .
\]
Choose in each $T_{n,j}$ a point $z_{n,j}$ such that 
\[
u_1 (z_{n,j}) \le \frac1{m (T_{n,j}) } 
\int_{T_{n,j}} u_1 (\zeta) dm (\zeta).
\]
It is enough to estimate this last average.  We apply Fubini's theorem:
\begin{multline*}
\frac1{m (T_{n,j}) } 
\int_{T_{n,j}} u_1 (z) dm (z)
=
\frac1{m (T_{n,j}) } 
\int_{T_{n,j}} 
\int_{\{w: d(z,w) \le \delta\}} \log \frac1{d(z,w)} d \nu(w)
dm (\zeta)
\\
\le \frac1{m (T_{n,j}) } 
\int_{Q_{n,j}} 
\left(
\int_{T_{n,j}}
\log \frac1{d(z,w)} dm (z)
\right)
d \nu(w)
\le
C  \int_{Q_{n,j}}  d \nu(w),
\end{multline*}
where the last inequality is due to the following explicit estimate :
for $w \in Q_{n,j}$, 
\begin{multline*}
\int_{T_{n,j}}
\log \frac1{d(z,w)} dm (z)
\le
\sum_{k=0}^\infty \log \frac1{2^{-k-1}} \,
m \left( \{ z : 2^{-k-1} < d(w,z) \le 2^{-k}\} \right) 
\\
\le
C \sum_{k=0}^\infty k 2^{-2n-2k} 
\le C 2^{-2n} \le C m (T_{n,j}).
\end{multline*}

Now we set 
\[
p_{n,j} := \frac1{|I_{n,j}|} \int_{ Q_{n,j}}  (1-|w|) d \nu(w) 
\approx \int_{ Q_{n,j}} d \nu(w).
\]
Since the ``squares" 
$Q_{n,j}$ are disjoint, the condition (\ref{Blaschkemeas}) implies that 
\begin{equation}
\label{sumatota}
\sum_{n,j} |I_{n,j}|p_{n,j} < \infty.
\end{equation}
Hence, $\{ p_{n,j} \}$ satisfies
(\ref{discNSC}) and Theorem \ref{discrete} provides a positive measure 
$\mu$ on the unit circle such that for any dyadic arc $I_{n,j}$, 
\begin{equation}
\label{discmaj}
\mu(I_{n,j}) \ge |I_{n,j}|p_{n,j} .
\end{equation}
However, whenever condition \eqref{sumatota} is verified, a more direct 
construction can be applied. Namely, one may take $d\mu = f(\eit) d 
\theta$ where
\[
f(\eit) := \sum_{I_{n,k} \ni \eit} p_{n,k} .
\]
Observe that $\int_0^{2\pi} f(\eit) d \theta = \sum |I_{n,k}| p_{n,k}$, 
and since $f(\eit) \ge p_{n,k}$ whenever $\eit \in I_{n,k} $, the measure 
$\mu$ satisfies \eqref{discmaj}.

Then, by the remarks 
before Theorem \ref{discrete}, there exists a harmonic function 
$h_1$ such that
\[
h_1 (z_{n,j}) \ge  C p_{n,j} \ge C \int_{ Q_{n,j}} d \nu(w)
\ge u_1 (z_{n,j}).
\]
\end{proof*}

\section{Proof of Propostion \ref{RnotLip}}
\label{PRNL}

The gist of this proof is that the Poisson kernel itself
cannot be log-Lipschitz with a constant better than $2$.

Denote
$f(z)=\Re\frac{1+z}{1-z} = P_z(1)$.
Fix  $\gamma >0$. For small $\delta>0$ consider the function
$$
g(z)=f((1-\delta)z),
$$
positive and harmonic in $\frac 1{1-\delta}\D$.
For $0<\varepsilon<\delta/2$ consider
$$
g_\varepsilon(z)=g\bigl(\frac {z}{1-\varepsilon}\bigr),\qquad z\in\D.
$$
Denote 
$$
M_\delta 
=\sup_{0<\varepsilon<\delta/2,\,z\in\D}\|\nabla(\log g_\varepsilon)(z)\|.
$$
Put
$$
h(z)=g_\varepsilon(z),\qquad z\in(1-\varepsilon)\T.
$$

{\bf Claim.} For sufficiently small $\varepsilon<\delta/2$,
for any $w\in (1-\varepsilon)\T$, $z\in\D$, then
$$
\left| \log g_\varepsilon(w)-\log g_\varepsilon(z) \right|
\le M_\delta |z-w|\le
\gamma \rho(z,w),
$$
where
$\rho$ is defined as in \eqref{Poincare}.

Applying this to the special case where $z \in (1-\eps)\T$, 
we see that for sufficiently small $\varepsilon$,
$\log h$ satisfies the Lipschitz condition with respect to
$\rho$ with constant $\gamma$.

Consider the $\gamma$-$\log$-Lipschitz (with respect to $\rho$) 
extension of $h$ to $\D$:
$$
H(z)=\sup_{w\in(1-\varepsilon)\T}h(w)e^{-\gamma\rho(w,z)}.
$$
It also 
follows from the Claim that for sufficiently small $\varepsilon<\delta/2$,
$$
H\le g_\varepsilon\qquad \text{on }\D.
$$

{\it Proof of the Claim.}

To prove the second inequality, note that if 
$\bigl|\frac{z-w}{1-z\bar w}\bigr|\le\frac 12$, then
$|1-z\bar w|\le u(\varepsilon)\to 0$ as $\varepsilon\to 0$,
and then (since $\log\frac{1+s}{1-s}\ge 2s$, $0\le s<1$)
$$
2\rho(w,z) 
\ge 2\bigl|\frac{z-w}{1-z\bar w}\bigr|\ge \frac 2{u(\varepsilon)}|z-w|. 
$$

If $\bigl|\frac{z-w}{1-z\bar w}\bigr|>\frac 12$, and
$|z-w|\le \frac{\gamma}{2M_\delta}$, then
$$
2\rho(w,z) \ge \log 3 > 1 \ge \frac {2M_\delta}{\gamma}|z-w|. 
$$

Finally, if $|z-w|>\frac{\gamma}{2M_\delta}$, then
$$
\bigl|\frac{z-w}{1-z\bar w}\bigr|\ge 1-
u_1\bigl(\varepsilon,\frac {\gamma}{2M_\delta}\bigr), 
$$
where $u_1\bigl(\varepsilon,\frac {\gamma}{2M_\delta}\bigr)\to 0$
as $\varepsilon\to 0$, and
$$
2\rho(w,z) \ge \log\frac{1}{u_1\bigl(\varepsilon,\frac {\gamma}{2M_\delta}\bigr)}
\ge \frac{2M_\delta}{\gamma}, 
$$
as $\varepsilon\to 0$. This completes the proof of our Claim.
\bigskip

By definition of $H$, we have $L_\gamma(H)=H$.

Since $R(H)$ is superharmonic, and $H$ coincides with $g_\varepsilon$
on $(1-\varepsilon)\T$, with $g_\varepsilon$ harmonic in $\D$, we obtain
$$
R(H)\ge g_\varepsilon\qquad\text{on }(1-\varepsilon)\D.
$$
Furthermore, since $H\le g_\varepsilon$ on $\D$,
$$
R(H)\le g_\varepsilon\qquad\text{on }\D.
$$
As a result,
$$
R(H)=g_\varepsilon\qquad\text{on }(1-\varepsilon)\D,
$$
and
$$
\frac{\partial}{\partial x}(\log R(H))(0)=
\frac{\partial}{\partial x}(\log g_\varepsilon)(0)=\frac
{1-\delta}{1-\varepsilon}\frac{\partial}{\partial x}(\log f)(0)=2\frac
{1-\delta}{1-\varepsilon}.
$$

Therefore, the function $R(H)$ is no better than 
$2\frac{1-\delta}{1-\varepsilon}$-$\log$-Lipschitz, and
$$
R(L_\gamma(H))=R(H)\ne R_c(H),\qquad c<2\frac{1-\delta}{1-\varepsilon}.
$$
This completes the proof.

\section{Comparison of upper envelopes}
\label{compenv}

The aim of this section is to prove Theorem \ref{twostepmaj}.  
Proposition \ref{liplarghmaj} shows that it will follow from the following 
result.

\begin{lemma}
\label{twostep}
Given $C_0>0$, there is a $C_1 \ge C_0$ such that if 
$\varphi$ is Log-Lipschitz with constant $C_0$, and admits
a non-trivial superharmonic majorant,
then there exists 
a superharmonic function $v$  
such that $v \ge \varphi$ and $v$ verifies \eqref{LipinLarg} with constant
$C_1$.
\end{lemma}

\begin{proof}
Let
$u \in Sp(\Dd)$ such that $u \ge \varphi$. We will notice 
(Lemma \ref{avlip}) that an averaged 
version of $u$ always satisfies (\ref{LipinLarg}), and that (up to a 
multiplicative constant) it provides the regular 
superharmonic majorant we are looking for.

For 
any $\delta \in (0, \frac12]$, let $D_H (z,\delta)$ 
stand, as above, for a hyperbolic disk of radius $\delta$ 
centered at $z$. Let $d\beta (z) := \diz^{-2} d m (z)$ be the invariant 
measure on the disk. For any M\"obius automorphism $\phi$ of the disk, 
any measurable function $f$ and any measurable set $E$, we have
(see e.g. \cite[(2.19), p. 19]{St})
\[
\int_E f \circ \phi \, d\beta = \int_{\phi(E)} f  \, d\beta.
\]
For any measurable function $g$ on the unit disk, let
\[
g_\delta (z) := \frac1{\beta (D_H (z,\delta))} \int_{D_H (z,\delta)} g \, d 
\beta.
\]

Since $\varphi$ is Log-Lipschitz with constant $C_0$, 
we have $\varphi (w) \ge e^{-C_0 \delta} \varphi 
(z)$ for any $w \in D_H (z,\delta)$, so that, being an average of such 
values, 
$\varphi_\delta (z) \ge e^{-C_0 \delta} \varphi (z)$. Now $\varphi \le u$ 
implies that 
\[
\varphi (z) \le  e^{C_0 \delta} \varphi_\delta (z) 
\le  e^{C_0 \delta} u_\delta (z) .
\]

The proof of Lemma \ref{twostep} will conclude with the next two Lemmas.
\end{proof}

\begin{lemma}
\label{avlip}
There exists an absolute constant $\kappa$ such that
for any positive valued superharmonic function $u$, and
for any $z$, $w$ in the unit disk such that $d (z,w) \le \delta/4$,
one has $u_\delta (w) \le \kappa u_\delta (z)$, and therefore 
$u_\delta$ verifies (\ref{LipinLarg}).
\end{lemma}

\begin{proof}
Recall that since $u$ is superharmonic, for any $z \in \mathbb D$ and 
$r_1<r_2$, 
\begin{equation}
\label{avdec}
\frac1{\beta(D_H (z,r_1) )} \int_{D_H (z,r_1)} u \, d \beta
\ge
\frac1{\beta(D_H (z,r_2) )} \int_{D_H (z,r_2)} u \, d \beta .
\end{equation}

Pick a constant $K < 1$ such that $\rho(z,w) \le \delta/4$ implies
that $D_H(w, K \delta ) \subset D_H(z, \delta )$.
Then \eqref{avdec} implies
\[
u_\delta (w) \le  u_{K_1 \delta} (w) \le
\frac1{\beta(D_H (w, K \delta))}  \int_{D_H (z, \delta)} u \, d 
\beta ,
\]
by the inclusion of discs and positivity of $u$. But
\[
\frac1{\beta(D_H (w, K \delta))}  \int_{D_H (z, \delta)} u \, d 
\beta \le
\kappa \frac1{\beta(D_H (z, \delta))}  \int_{D_H (z, \delta)} u \, d 
\beta = \kappa u_\delta (z),
\]
since $\beta(D_H (w,K \delta))$ and $\beta(D_H (z,\delta))$ are comparable.
\end{proof}

\begin{lemma}
\label{avsupharc}
Let $u$ be a positive superharmonic function on the disc. Then $u_\delta$ 
is also superharmonic.
\end{lemma}

\begin{proof}
We will follow the 
notations and use the results of \cite[Chapter 4, pp. 34--39]{St}, itself 
inspired by \cite{Ul}. Let $\phi_z$ denote the unique involutive M\"obius 
automorphism of the disk which exchanges $z$ and $0$.
Define the invariant convolution of two measurable 
functions by
\[
f * g (z) := \int_{\Dd} f(w) g\circ \phi_z (w) \, d\beta (w),
\]
whenever the integral makes sense. This operation is commutative 
\cite[bottom of page 34]{St}.
If we set
\[
\Omega_\delta (z) := \frac1{\beta(D_H (0,\delta))} \chi_{D_H (0,\delta)} 
(z),
\]
then $f_\delta = f * \Omega_\delta$. Note that $D_H 
(0,\delta)=D(0,\tanh \delta)$ so that $\Omega_\delta$ is a radial 
function.

The $\mathcal M$-subharmonic functions defined in \cite[Chapter 4, (4.1)]{St} 
reduce for $n=1$ to ordinary subharmonic functions.
The invariant Laplacian (Laplace-Beltrami operator for the Bergman metric 
of the ball) reduces in the case $n=1$ to
\[
\tilde \Delta = 2 \diz^2 \frac{\partial^2}{\partial z \partial \bar z},
\]
so that $\mathcal C^2$ superharmonic functions $g$ can be characterized as those 
such that $\tilde \Delta g \le 0$.  Since our function $\Omega_\delta$ is 
not smooth, we need to perform an approximation argument.
It will be enough to show that $u * \Omega_\delta$ can be approximated 
from below by an increasing sequence of
$\mathcal C^2$ superharmonic functions. Pick an
increasing sequence of 
smooth, nonnegative, {\em radial} functions $\Omega_{\delta,n}$ so that 
$\lim_{n\to\infty} \Omega_{\delta,n} = \Omega_{\delta}$ almost everywhere. 
Then the monotone convergence theorem tells us that $u*\Omega_{\delta,n}$ 
converges to
$u * \Omega_\delta$, and the sequence is clearly increasing. For $f \in 
\mathcal C^2 (\Dd)$, by \cite[(4.11), p. 36]{St}, 
\[
(\tilde \Delta f) (a) =
\lim_{r\to 0} \frac4{r^2} [ (f*\Omega_r)(a) - f(a) ].
\]
Now,  twice applying
Ulrich's lemma about associativity of the invariant convolution 
when the middle element is radial \cite[Lemma 4.5, p. 36]{St}, we have
\begin{multline*}
\tilde \Delta (u*\Omega_{\delta,n}) =
\lim_{r\to 0} \frac4{r^2} \left[ (u*\Omega_{\delta,n})*\Omega_r -
u*\Omega_{\delta,n} \right]
= \lim_{r\to 0} \frac4{r^2} \left[ u* (\Omega_{\delta,n}*\Omega_r) -
u*\Omega_{\delta,n} \right]
\\
= \lim_{r\to 0} \frac4{r^2} \left[ u* (\Omega_r*\Omega_{\delta,n}) -
u*\Omega_{\delta,n} \right]
=  \lim_{r\to 0} \frac4{r^2} \left[ (u* \Omega_r)*\Omega_{\delta,n} -
u*\Omega_{\delta,n} \right]
\\
=  \lim_{r\to 0} \frac4{r^2} \left[ (u* \Omega_r - u)*\Omega_{\delta,n} 
\right].
\end{multline*}
Since $u$ is superharmonic, $u* \Omega_r - u \le 0$, and since 
$\Omega_{\delta,n} \ge 0$, we finally have 
$\tilde \Delta (u*\Omega_{\delta,n}) \le 0$.
\end{proof}

\section{Counterexamples}

We recall that for $z$, $w \in U_+$, the Gleason distance is given by 
$d (z,w) = \left| \frac{z-w}{z-\bar w} \right|$, and again
$\rho= \frac12 \log \frac{1+d}{1-d}$. For a sequence $\{z_k = 
x_k+iy_k\}$, the Blaschke condition reads 
\[
\sum_k \frac{y_k}{1+|z_k|^2} < \infty,
\]
which reduces to $\sum_k y_k < \infty$ when the sequence is bounded.
Finally, the Blaschke product with zeroes at the $z_k$ is
\[
B(z) := \prod_k \frac{z-z_k}{z-\overline z_k}.
\]

\subsection{Construction of the examples}

We shall use a monotone increasing function $f$ from $(0,1]$ to $(0,1]$ 
satisfying $f(t)\le t$, $\lim_{t\to 0} f(t)/t =0$, and there exists a constant
$C_0 >1$ such that $f(t/2) \ge f(t)/C_0$ and $f(3t/2) \le C_0 f(t)$. 
For instance, any 
function of the form $f(t) = t^m$, for $m>1$ will satisfy those conditions.

Define an arc $\gamma$ in $U_+$ by $\gamma(x) := x+if(x)$, $0<x\le1$. 
Pick a sequence of points $a_k := x_k + i f(x_k)$ on $\gamma$ such that 
$x_0=1$, $\{x_k\}$ decreases towards $0$, and $\rho (x_k,x_{k+1}) = 2 
\delta>0$ is constant, with $\delta <1/2$ to be specified later. 

It is easy to see that the sequence $\{a_k\}$ is interpolating (each 
Carleson window of size $f(x_k)$ around the point $a_k$ can contain at 
most a fixed finite number of other points in the sequence), so the 
Blaschke product $B$ associated to the sequence will verify
\begin{equation}
\label{BlasEst}
\log \frac{1}{|B(z)|} \le C_\delta + \log \frac1{d(a_k,z)} ,
\mbox{ for any } z \in D_H (a_k, \delta).
\end{equation}

Now, given a decreasing function $\sigma$ from $\Rr_+^*$ to $\Rr_+^*$ 
with $\lim_{t\to 0} \sigma(t) = \infty$, if 
we set $\varphi_\sigma (x+iy) := \sigma (|x|+y)$, we have 
$M\varphi_\sigma (\xi) = \sigma(|\xi|)$, therefore
$|\{  M \varphi > t \}| = 2 \sigma^{-1} (t)$ (where $\sigma^{-1}$ stands 
for the inverse function).

\begin{lemma}
\label{constsuperh}
There exists $k_0 \in \Zz_+$ such that if $B_0$ stands for the Blaschke 
product associated to the sequence $\{a_k, k\ge k_0\}$, then
the following function
\begin{eqnarray*}
\varphi(z) & := & \log \frac1{|B_0(z)|}, \mbox{ if } z \notin
\bigcup_{k \ge k_0} D_H (a_k, \delta) ,\\
\varphi(z) & := & \min \left( \log \frac1{|B_0(z)|} , \min_{D_H (a_k, 
\delta/2)} \varphi_\sigma \right), \mbox{ if } z \in D_H (a_k, \delta) 
\mbox{ for some } k \ge k_0 
\end{eqnarray*}
is superharmonic on $U_+$.
\end{lemma}

Notice that the definition of $\varphi$ and the fact that $ \log 
\frac1{|B_0|}$ is bounded outside of 
$\bigcup_{k \ge k_0} D_H (a_k, \delta)$
imply that $\varphi \le \max (C_1, \varphi_\sigma)$ for some $C_1 >0$, and 
therefore for $t > C_1$,
\begin{equation}
\label{maxfcnest}
\left| \{ \xi \in \Rr : M \varphi (\xi) > t \} \right| \le 2 \sigma^{-1} (t) .
\end{equation}

\begin{proof*}{\bf Proof of Lemma \ref{constsuperh}:}

Since $ \log \frac1{|B_0|}$ is superharmonic, it will be enough to show 
that 
\begin{equation}
\label{ineqBphi}
\min \left( \log \frac1{|B_0(z)|} , \min_{D_H (a_k, 
\delta/2)} \varphi_\sigma \right) = \log \frac1{|B_0(z)|}
\end{equation}
 in a neighborhood 
of $\partial D_H (a_k, \delta)$ for each $k \ge k_0$. For 
$\rho(z,a_k) \ge (3/4) \delta$, we have from (\ref{BlasEst}):
\[
 \log \frac1{|B_0(z)|} \le C_\delta + \log \frac1{\tanh \frac34 \delta},
\]
while $\lim_{k\to\infty} \min_{D_H (a_k, \delta/2)} \varphi_\sigma = \infty$
because $\lim_{t\to 0} \sigma(t) = \infty$. So there exists $k_0$ so that 
(\ref{ineqBphi}) is achieved for $3\delta /4 \le d(a_k,z)\le \delta$.
\end{proof*}

\subsection{Proofs}
 \ {}

\begin{proof*}{\bf Proof of Proposition \ref{sharpmaxf}:}

Suppose there exists some $h \in Ha(U_+)$ 
such that $h(z) \ge \varphi(z) \ge 0$, for any $z \in U_+$. 
Then Harnack's inequality implies that if $z$, $w \in U_+$ are such that
$\rho(z,w) < 2 \delta$, we have $h(w) \ge C'_\delta h(z)$.

Provided $k_0$ is large enough, we have $|x|+y \le 2 x_k$ for any
$x+iy \in D_H (a_k, \delta/2)$, so that 
$\min_{D_H (a_k, \delta/2)} \varphi_\sigma
\ge \sigma (2x_k)$. Harnack's inequality as used above then yields 
\[
h(z) \ge C'_\delta \sigma (2x_k), \mbox{ for all }
z \in D_H (a_k, 2\delta) .
\]
Take any $x \le x_{k_0}$. There is a unique $k \ge k_0 +1$ such that
$x_{k+1} \le x < x_k$. 
Since the whole curve $\gamma$ is covered by the disks $D_H (a_j, 
2\delta)$ we have either 
$x \in D_H (a_{k+1}, 2\delta)$ and 
$h(z) \ge C'_\delta \sigma (2x_{k+1}) 
\ge C'_\delta \sigma (2x)$, or $x \in D_H (a_{k}, 2\delta)$ and 
$h(z) \ge C'_\delta \sigma (2x_{k})$ ; in this latter case we have 
\[
x_k > x \ge x_k - C f(x_k) \ge \frac12 x_k
\]
for $k_0$ large enough. So in each case we have 
\begin{equation}
\label{hmesgran}
\frac1{C'_\delta} h(x+if(x)) \ge \sigma (4x), \mbox{ for any }x \le x_{k_0}.
\end{equation}

On the other hand, there exists $c>0$ and a positive measure $\mu$ on the 
real line with $\int \frac{d\mu(t)}{1+t^2} < \infty$ such that
\[
h_1(x+iy) :=
\frac1{C'_\delta} h(x+iy) = cy + \int_{\Rr} P_y (x-t) d\mu (t),
\]
where $P_y (t):= \frac1\pi \frac{y}{y^2+t^2}$ stands for the Poisson 
kernel. In particular, if $|x| \le 1$, $y\le 1$, 
\[
h_1(x+iy) \le \int_{-2}^2 P_y (x-t) d\mu (t) + c + 
\int_{|t|\ge 2} \frac1{1+t^2} d\mu (t) 
\le \int_{-2}^2 P_y (x-t) d\mu (t) + C_2 .
\]
Therefore, again for $k_0$ large enough, (\ref{hmesgran}) implies
\[
\int_{-2}^2 P_{f(x)} (x-t) d\mu (t) \ge h_1(x+if(x)) - C_2 
\ge  \sigma (4x) -C_2 \ge  \frac12 \sigma (4x).
\]
Hence
\begin{equation}
\label{intsigbdd}
\int_0^1 \sigma (4x) dx 
\le
2 \int_{-2}^2 \left( \int_0^1 P_{f(x)} (x-t) dx \right) d\mu (t).
\end{equation}

Now we choose $f(x)=x^2$ (the choice of $f$ is inessential here). 
We claim that $\int_0^1 P_{x^2} (x-t) dx \le C$. Indeed,
let 
\begin{eqnarray*}
\int_0^1 &&\frac{x^2}{ (x-t)^2+x^4} dx \\
&=&
\int_{0 \le x \le 1, |x-t| \le t^2/4}
+
\int_{0 \le x \le 1, t^2/4 \le |x-t| \le t}
+
\int_{0 \le x \le 1, t \le |x-t|}
\frac{x^2}{ (x-t)^2+x^4} dx \\
&=&: I + II + III.
\end{eqnarray*}
Then
\[
I \le \int_{ |x-t| \le t^2/4} \frac{dx}{x^2}
\le 2 (t^2/4) (4/t^2) \le 2,
\]
\[
II \le \int_{t^2/4 \le |x-t| \le t} \frac{x^2dx}{(x-t)^2} 
\le Ct^2 \int_{t^2/4 \le |x-t| } \frac{dx}{(x-t)^2} 
\le C t^2 ( t^2/4 )^{-1} < \infty,
\]
and, since $t \le |x-t|$ implies that $|x-t| \ge x/2$, 
\[
III \le \int_{0 \le x \le 1}
\frac{x^2}{ (x/2)^2+x^4} dx  \le 4.
\]
The claim is proved.

Finally, (\ref{intsigbdd}) together with the claim says that
\[
\int_0^1 \sigma (4x) dx  \le C \int_{-2}^2 d\mu (t) < \infty,
\]
therefore $\int_0^1 \sigma (x) dx < \infty$. 

Now we choose $\sigma := s^{-1}$, where $s$ is the function given in 
Proposition \ref{sharpmaxf}. Since $\int_0^1 \sigma (x) dx \ge 
\int_{s^{-1}(1)}^\infty s(t) dt$, we have a contradiction with 
$\int_1^\infty s(t) dt = \infty$.
\end{proof*}

\begin{proof*}{\bf Proof of Proposition \ref{anyrate}:}

We construct a function $\varphi$ as above with $\sigma (t) = t^{-2}$.
If there was a harmonic $h$ such that $h \ge \varphi$, Harnack's inequality 
would imply as before (\ref{hmesgran}), which is impossible because we would have $Mh(\xi) \ge C 
\xi^{-2}$ for $\xi >0$ and small enough, which is not a weakly 
integrable function. 

We must now check the growth property for $\varphi$. 
As we have seen before, if $z \notin
\bigcup_{k \ge k_0} D_H (a_k, \delta)$, then $\varphi(z)$ is bounded by a 
constant. So we can restrict ourselves to points in a hyperbolic 
neighborhood of the curve $\gamma$.

For points $z \in \gamma$ with $z=x+if(x)$, $x\le x_{k_0}$, 
\[
\varphi (z) \le \varphi_\sigma (x+if(x)) = \sigma (x+f(x)) \le \sigma(x),
\]
therefore on the curve $\gamma$ we have $\varphi (x+iy) \le \sigma 
(f^{-1}(y))$.

In general, there is a $x'$ such that $z = x+iy \in D_H (x'+if(x'), \delta)$,
and $\varphi(z) \le \varphi(x'+if(x')) \le \sigma(x')$.
Now $y \le C f(x')$, so
\[
\varphi(z) \le \sigma(x') \le \sigma (f^{-1}(y/C)) .
\]
With our choice of $\sigma$, one can check that choosing any $f$ so that
\[
f(\xi) \le \frac1C s^{-1}(\frac1{\xi^2})
\]
will yield $\varphi(z) \le s(y)$, q.e.d.
\end{proof*}

\vskip1cm
Alexander Borichev, LaBAG, Universit\'e Bordeaux 1, 351 cours de la Lib\'eration
33405 Talence CEDEX, France.

Alexander.Borichev@math.u-bordeaux.fr

\vskip.3cm

Artur Nicolau, Departament de Matem\`atiques, Universitat Aut\`onoma
de Bar\-ce\-lo\-na, 08193 Bellaterra, Spain.

artur@mat.uab.es
\vskip.3cm

Pascal J. Thomas, Laboratoire de Math\'ematiques Emile Picard, UMR CNRS 5580,
Universit\'e Paul Sabatier, 118 route de Narbonne, 31062 Toulouse CEDEX, France.

pthomas@cict.fr

\end{document}